\documentclass[a4paper,12pt]{article}
\usepackage{amsmath}
\usepackage{amssymb,xypic,geometry,latexsym}
\usepackage{citesort}

\usepackage{xypic}
\author{ A. Mutlu and T.Porter}
\date{}

\title{Crossed squares and 2-crossed modules}

\newtheorem{defn}{Definition}[section]
\newtheorem{prop}[defn]{Proposition}
\newtheorem{thm}[defn]{Theorem}

 \newenvironment{pf}{{\bf Proof:}}{$\Box$\mbox{}}

\begin{document}

\maketitle
\begin{abstract}

{\noindent\bf A. M. S. Classification:}  18G30, 18G55. 
\end{abstract}


\section{Introduction}

Simplicial groups were first studied by D. M. Kan in the 1950s \cite{kan1}. Early
work by Kan himself, Moore, Milnor, and Dold showed that\\
(a) these objects have a well structured homotopy theory, \\
(b) they modelled all homotopy types of connected spaces,\\
(c) abelian simplicial groups were an equivalent tool to that of chain complexes
and could therefore be applied within homological algebra, and\\
(d) in low dimensions, calculations were possible, provided, for instance, the 
simplicial group was free with chosen $`CW$-basis'. 

Simplicially enriched groupoids have a more recent birth, but have the same
sort of attributes plus being able to model non-connected homotopy types.  The 
shortened form of their name `simplicial groupoid' used by Dwyer and Kan
\cite{D&K} is more usually used for these objects although not strictly
correct as simplicial objects in the category of groupoids form a much larger
setting than do simplicially enriched groupoids.  None the less we will often
use that shortened form here.

Crossed modules model homotopy 2-types. Crossed squares model 3-types. Crossed 
$n$-cubes model $(n+1)-$types, cf. \cite{porter} and the references therein.
Conduch\'e \cite{conduche} has an alternative model for 3-types namely
2-crossed modules and Baues \cite{baues}  uses a variant of this, the
quadratic module, in some of his work.  Another model for 3-types was
introduced by Brown and Gilbert, \cite{browngilbert}, (involving a particular
subdivided triangular diagram in its derivation). This model, braided
crossed modules and their lax counterpart, the braided categorical groups, have
been further studied  by members of the Granada Algebra group (cf. Carrasco and
Cegarra, \cite{C&C2}, and Garzon and Miranda, \cite{garzonmiranda}, for
example).  Similar object, Gray groupoids, have been studied by Joyal and Tierney
(unpublished) and have recently started to receive some attention in the TQFT literature.
The proof of the equivalence between Gray groupoids (that is, 2-groupoid enriched groupoids) and Conduch\'e's models is discussed in \cite{K&P} and more detailed references to other work by LeRoy, Berger, Marty, and, of course, the original proof from the 1980s by Joyal and Tierney, can be found there as well as in the bibliography of this paper.

In a letter to Brown and Loday, dated in the mid 1980s, Conduch\'e pointed out 
that given a crossed square,
$${\cal M} = \left(\diagram
L\ar[r]\ar[d] &M\ar[d]\\
                             N\ar[r] & P
\enddiagram
\right),$$
the mapping cone complex of $\cal M$,
$$L\rightarrow M\rtimes N\rightarrow P$$ constructed by Loday in \cite{loday}, 
has a 2-crossed module structure.  He ended his letter by pointing out that
this seemed to give a canonical and direct way to link the categories of crossed squares 
and 2-crossed modules, and his results suggested that it should yield some
sort of equivalence between the two models.

Although much more is known on this area than when Conduch\'e's letter was
written, the background underlying structure still seems obscure.  In this
paper we aim to shed some light on the 2-crossed module structure given by
Conduch\'e and also on the subdivided triangular diagram of Brown and
Gilbert.  This, in fact, provides the key and suggests ways of generalising
Conduch\'e's construction to higher $n$-types.

\section{Preliminaries}
\subsection{Simplicial Groups and Groupoids}
We assume that the reader is conversant with the basic theory of simplicial
sets and simplicial groups.  The following merely sets up notation and some
conventions.

Let  $\mathfrak{Grp}$ be the category of groups. A {\em simplicial group}
${\bf G}$ consists of a family of  groups $\{G_n\}$ together with face and degeneracy maps
$d_i=d_i^n: G_n\longrightarrow G_{n-1},$ ~$0\leq i< n$ ~$(n\ne 0)$ and
$s_i=s_i^n:G_n\longrightarrow G_{n+1}$, $0\leq i\leq n,$ satisfying the usual
simplicial identities given in \cite{curtis}, \cite{kan1} and \cite{kan2}. It
can be completely described as a functor ${\bf G}: \Delta^{op}\longrightarrow
\mathfrak{Grp}$ where $\Delta$ is the category of finite ordinals, ${  [n]}=\{0<1<2<\cdots <n\}$, and increasing maps. We will denote the category
of simplicial group by $\mathfrak{SimpGrps}.$ We have for each $k\geq 0$, a
subcategory $\Delta_{\leq k}$ determined by objects $[j]$ of $\Delta$ with
$j \leq k.$ A $k$-truncated simplicial group is a functor from $(\Delta_{\leq
  k}^{op})$ to $\mathfrak{Grp}$, where $(\Delta_{\leq
  k}^{op})$ has the obvious meaning.

\medskip

\textbf{Remark:}

We will restrict detailed attention in the main to simplicial groups
and hence to connected homotopy types.  This is traditional but a bit
unnatural as all the results and definitions extend with little or no
trouble to simplicial groupoids in the sense of Dwyer and Kan \cite{D&K} and
hence to non-connected homotopy types.  It should be noted that such
simplicial groupoids have a fixed and constant simplicial set of objects and
so are not merely simplicial objects in the category of groupoids.  In this
context if $\mathbf{G}$ is a simplicial groupoid with set of objects $O$, the
natural form of the Moore complex $\mathbf{NG}$ (see below) is given by the same formula
as in the reduced case, interpreting Ker$ d^n_i$ as being the subgroupoid of
elements in $G_n$ whose $i$\textsuperscript{th} face is an identity of $G_{n-1}$.
Of course if $n \geq 1$, the resulting $NG_n$ is a disjoint union of groups,
so $\mathbf{NG}$  is a disjoint union of the Moore complexes of the vertex
simplicial groups of $\mathbf{G}$ together with the groupoid $G_0$ providing
elements that allow conjugation between (some of) these vertex complexes
(cf. Ehlers and Porter \cite{ep}).

\medskip

Consider the product $\Delta\times\Delta$ whose objects are pairs $({ [p]},{
  [q]})$ and whose maps are pairs of weakly increasing maps. A functor
$\mathcal{G}:(\Delta\times\Delta)^{op}\longrightarrow\mathfrak{Grp}$ is called
a {\em bisimplicial group}. To give $\mathcal{G}$ is equivalent to giving for
each $(p,q)$ a group $G_{p,q}$ and homomorphisms
$$
\begin{array}{llllc}
d_i^h&:&G_{p,q}\longrightarrow G_{{p-1},q}&\\
s_i^h&:&G_{p,q}\longrightarrow G_{{p+1},q}&\qquad i:0,1,\cdots,p\\
d_j^v&:&G_{p,q}\longrightarrow G_{p,{q-1}}&\\
s_j^v&:&G_{p,q}\longrightarrow G_{p,{q+1}}&\qquad j:0,1,\cdots,p\\
\end{array}
$$
such that the maps $d_i^h,~s_i^h$ commute with $d_j^v, ~s_j^v$ and $d_i^h,
~s_i^h$ (resp. $d_j^v,~ s_j^v$) and  satisfy the usual simplicial identities. Here
$d_i^h,~s_i^h$ denote the horizontal operators and $d_j^v~s_j^v$  denote the
vertical operators. A bisimplicial group can also be thought of as a
simplicial object in the category of simplicial groups.

Multisimplicial groups are similarly defined. Explicitly if $n$ is a positive
integer, form the $n$-fold product $\Delta^{\times n}$ of $\Delta$ with
itself, then an $n$-simplicial group is merely a functor
$\mathcal{G}:(\Delta^{\times n})^{op}\longrightarrow\mathfrak{Grp}$. Of course 
an $n$-simplicial group leads to an $n$-indexed family of
groups $G_{p_1,\ldots, p_n}$ with face and degeneracy operators in the obvious 
way, generalising the case $n = 2$ discussed above.

All of this extends painlessly to simplicial groupoids and with obvious
definitions of morphisms of $n$-simplicial groups and groupoids, we thus get a
whole host of categories  $n-\mathfrak{SimpGrps}$ and $n-\mathfrak{SimpGrds}$.

Recall from \cite{porter} that a {\em normal chain complex} of groups,
$(X,d)$ means one in which each Im$~d_{i+1}$ is a normal subgroup of the
corresponding $X_i.$ Given any normal chain complex $(X,d)$ of groups and
an integer $n$, the truncation, $t_{n]}X$ of $X$ at level $n$ is defined by
\begin{displaymath}
(t_{n]}X)_i=\left\{\begin{array}{lll}
X_i & \text{if $i<n$}\\
X_i/\text{Im}d_{n+1}& \text{if $i=n$}\\
0 &\text{if $i>n.$}
\end{array}\right.
\end{displaymath}
The differential $d$ of $t_{n]}X$ is that of $X$ for $i<n,$ whilst $d_n$ is induced
by the $n^\mathrm{th}$ differential of $X$ and all other differentials are zero.

There are multicomplex and groupoid versions of this.  The only points worth
noting here are that any non trivial normal subgroupoid is a kernel and hence
must be a disjoint union of its vertex groups and that in a normal
$n$-multicomplex of groups, $(X, d_1, \ldots , d_n)$, each image in each
direction is assumed normal in each dimension.

\subsection{The Moore complex of simplicial group}

The archetypal normal $n$-complex  of group(oid)s is the Moore complex of an
$n$-simplicial group(oid). In particular given a simplicial group ${\bf G},$
{\em the Moore complex} $({\bf NG},\partial)$ of ${\bf G}$ is the normal chain complex defined by 
\[{\bf NG}_n =\bigcap_{i=0}^{n-1}\mbox{Ker}d_i \] and with differential $\partial_n : NG_n 
\to NG_{n-1}$ induced from $d_n$ by restriction. 

To form the Moore $n$-complex of an $n$-simplicial group(oid) one just
applies the functor $N$ repeatedly in each direction in turn.

There is an extensive theory of the Moore complex and its links with homotopy
theory.  We mention that $\pi_n(G)$ can be calculated by calculating the
$n^\mathrm{th}$ homology group of $({\bf NG},\partial)$. Explicitly
the {\em $n^\mathrm{th}$ homotopy group} $\pi_n({\bf G})$ of ${\bf G}$ is
defined to be the $n^\mathrm{th}$ homology of the Moore complex of ${\bf G},$ i.e.,
\[\pi_n({\bf G})\cong H_n({\bf NG},\partial)=\bigcap_{i=0}^{n}\mbox{Ker}d_i^n/d_{n+1}^{n+1}(\bigcap_{i=0}^{n}\mbox{Ker}d_i^{n+1}).\]
We say that the Moore complex ${\bf NG}$ of a simplicial group is of {\em
  length} $k$ if $NG_n=1$ for all $n\geq k+1$ so that a Moore complex of
length $k$ also has length $l$ for $l\geq k.$

A simplicial map $f:{\bf G}\to{\bf {G'}}$ is called an $n$-{\em equivalence} if it induces isomorphisms
\[\pi_k({\bf G})\cong\pi_k({\bf G'})\qquad \text{for}\quad k\leq n.\]
Composites of $n$-equivalences are, of course, also $n$-equivalences.
Two simplicial groups, ${\bf G}$ and ${\bf {G'}}$,  are said to have the same
$n$-{\em type} if there is a zig-zag chain of $n$-equivalences linking them. A simplicial
group ${\bf G}$ \emph{is an $n$-type} if $\pi_i({\bf G})=1$ for $i>n.$

The Moore complex carries a lot of fine structure and this has been studied,
e.g. by Carrasco and Cegarra, \cite{C&C1}, Wu \cite{wu}, and the present
authors in earlier papers in this series, \cite{mp1,mp2,mp3,mp4,mp5}.  The specific
structure of the $k$-truncation of the Moore complex for $k=1$ and 2 is now
well known.  For  $k=1$ this gives a crossed module, for  $k=2$, a 2-crossed
module, cf. Conduch\'e, \cite{conduche}.  We summarise his theory below.  The
hypercrossed complex structure of  Carrasco and Cegarra, \cite{C&C1}, clearly
models all homotopy types and on $k$-truncation, all $(k + 1)$-types, but the
structure does get unwealdy for $k$ larger than about 3 or 4.

All of this theory works for both simplicial groups and simplicially enriched 
groupoids with virtually identical presentations. Rather than write
`simplicial group(oid)' all the time we have written the exposition in terms
of simplicial groups but the other case works in the same way.
\section{2-crossed modules and simplicial groups} 
Crossed module techniques give a very efficient way of handling information
about a homotopy type. They correspond to 2-types (see \cite{conduche} and \cite{mp2}). As mentioned  above Conduch\'{e},
\cite{conduche}, in 1984 introduced the notion of 2-crossed module as a model for 3-types.

Throughout this paper we denote an action of $p\in P$ on $m\in  M$ by $p\cdot m={}^{p}m.$

A {\em crossed module} is a group homomorphism $\partial: M\to P$ together
with an action of $P$ on $M$ satisfying \\
(C1)\quad$\partial({}^{p}m)=p\partial(m)p^{-1}$ \\and\\ (C2)\quad
${}^{\partial{m}}{m'}=m{m'}m^{-1}$ for all $m,{m'}\in M ,~p\in P.$ 

\noindent This second condition is called the {\em Peiffer identity}. We will denote such a crossed
module by $(M,P,\partial).$ 

A {\em morphism} of crossed modules from $(M,
P,\partial)$ to $({M'},{P'},{\partial'})$ is a pair of group homomorphisms, $\phi:
M\to {M'}$, $ \psi:P\to {P'}$ such that $\phi({}^{p}m)={}^{\psi(p)}\phi(m)$ and ${\partial'}\phi(m)=\psi\partial(m)$. We thus get a category $\mathfrak{XMod}$ of crossed modules.

\medskip

{\em Examples of crossed modules}

(a)\quad Any normal subgroup $N$ in $P$ gives an inclusion map, inc$:N\to P$ which is a crossed module. Conversely given any arbitrary crossed module $\partial: M\to P,$ one can easily see that the Peiffer identity implies that $\partial M=P$ is a normal subgroup in $P.$\\
(b)\quad Given any $P$-module, $L$, the trivial homomorphism ${\bf 1}: L\to P$
is a crossed $P$-module for the given action of $P$ on $L$.

The following definition of 2-crossed modules  is equivalent to that given by
D.Conduch\'e, \cite{conduche}.

\medskip

\textbf{Definition:}

A \emph{2-crossed module}  consists of a complex of groups 
$$
\diagram
L\rto^{\partial_2} & M\rto^{\partial_1}& N 
\enddiagram
$$
together with an action of $N$ on $L$ and $M$ so that 
$\partial _2,\partial _1$ are morphisms of $N$-groups, where the group 
$N$ acts on itself by conjugation, and an $N$-equivariant function 
$$
\{\quad , \quad \}: M \times M \to L \\
$$
called a \emph{ Peiffer lifting}, which satisfies the following axioms:
$$
\begin{array}{llrcl}
2CM1: &  & \partial _2\{m, {m^{\prime}}\} & = & ({~}^{\partial_1
m}{m'}){~} (m({m'})^{-1}m^{-1}), \\ 
2CM2: &  & \{\partial _2(l), \partial _2({l^{\prime}})\} & = & [{l'},
l], \\ 
2CM3: & (i) & \{mm', {m''}\} & = & {~}^{\partial_1 m}
\{{m'}, {m''}\} \{m,{m'}{m''}({m'})^{-1}\}, \\ 
&(ii) & \{m, {m'}{m''}\} & = &\{m, {m'}\}
{~}{}^{mm'(m)^{-1}} \{m,{m''} \}, \\
2CM4: &  & \{m,\partial_2{l} \}\{\partial_2l,m\} & = & 
 {~}^{\partial_1 m}l(l)^{-1},\\
2CM5& & {~}^{n}\{m, {m^{\prime}}\} & = & \{{~}^{n}m, {~}^{n}{m^{\prime}}\},
\end{array}
$$
for all $l,{l^{\prime}}\in L,~m,{m^{\prime}}, {m^{\prime\prime}}
\in M$ and $n\in N.$

Here we have used ${}^ml$ as a shorthand for $\{\partial_2l,m\}l$ in
condition 2CM3(ii) where $l$ is $\{m, m^{\prime\prime}\}$ and $m$ is
  $mm'(m)^{-1}$. This gives a new action of $M$ on $L$.  Using this notation,
  we can split 2CM4 into two pieces, the first of which is tautologous:
$$\begin{array}{llrcl}
2CM4: & (a) & \{\partial_2l, m\} & = &{}^{m}(l). l^{-1}, \\
& (b) & \{m,  \partial_2l\} & = & ({}^{\partial_1m}l)({}^{m}l^{-1}).
\end{array}$$
The old action of $M$ on $L$, via $\partial_1$ and the $N-$action on $L$, is in general distinct 
from this second action with $\{m,  \partial _2(l)\}$ measuring the difference 
(by 2CM4(b)).  An easy argument using 2CM2 and 2CM4(b) shows that with this
action, ${}^ml$, of $M$ on $L$, $(L,M,\partial_2)$ becomes a crossed module.

We denote such a 2-crossed module  by $\{L,\,M,\,N,\,%
\partial _2,\,\partial _1\}.$
A morphism of 2-crossed modules is given by a  diagram
$$
\diagram
L \rto^{\partial_2} \dto_{f_2} & M \dto_{f_1} \rto^{\partial_1} & N 
\dto_{f_0} \\ {L'} \rto_{\partial'_2} & {M'}
\rto_{\partial'_1}& {N'}
\enddiagram
$$
where $f_0\partial _1=\partial _1^{\prime }f_1,\ 
f_1\partial _2=\partial_2^{\prime }f_2,$  
$$
f_1({}^{n} m_1)={}^{f_0(n)}f_1(m_1),\quad f_2({}^{n}
l)={}^{f_0(n)} f_2(l), 
$$
and 
$$
\{\quad , \quad \}f_1\times f_1=f_2\{\quad , \quad \}, 
$$
for all $l\in L,\ m_1\in M,\ n\in N.\ $ These compose in an obvious way. 

The groupoid analogues of these definitions are left to the reader. We will
concentrate on the reduced case i.e. with groups rather than groupoids. 

We thus can consider the category of 2-crossed modules denoting it as 
$\mathfrak{X_{2}Mod}.$ Conduch\'{e} \cite{conduche} proved that 
2-crossed modules give algebraic models of connected homotopy 3-types.
  \begin{thm} {\rm{(\cite{conduche}, \cite{mp3}})}
The category, $\mathfrak{X_2Mod}$, of $2$-crossed modules is equivalent to the category $\mathfrak{SimpGrp}_{\leq 2}$ of
simplicial groups with Moore complex of length $2$.\hfill $\square$
\end{thm}

\section{Cat$^2$-groups and crossed squares}
The following definition is due to D. Guin-Walery and J.-L. Loday, see  \cite{wl} and also \cite{loday}.

\textbf{Definition:}

{\rm A crossed square of groups} is a commutative square of groups 
$$
\diagram
L \rto^{\lambda} \dto_{\lambda'}  & M\dto^{\mu} \\
{N} \rto_{\mu'} & P  \\
\enddiagram
$$            
together with actions of $P$ on $L$, $M$ and~${N}.$ There are thus 
actions of ${N}$ on ~$L$ and ~$M$ via~ ${\mu'}$ 
and $M$ acts on $L$  and ${N}$ via $\mu$ and a function 
$h:M \times {N} \longrightarrow L$ such that,
for all $l \in L,~ m_1, m \in M, ~{n_1},{n} \in {N}$ ~and~ $p \in
P$  the following axioms hold:
\begin{enumerate}
\item  the homomorphisms $\lambda ,\ \lambda ^{\prime },\ \mu \ \mu'$ and 
$\kappa = \mu \lambda  =\mu' {\lambda}^{\prime }$ are crossed modules for the
corresponding actions and 
the morphisms of maps $(\lambda) \to (\kappa); ~(\kappa) \to(\mu)$;~
$(\lambda') \to (\kappa)$; ~and ~ $(\kappa) \to(\mu')$
are morphisms of crossed modules,

\item  $\lambda h(m,n)=m{}{}^{\mu'({n})}m,$

\item  ${\lambda'}h(m, n)= {}^{\mu(m)}n(n)^{-1},$

\item  $h(\lambda(l),n)=l~{}^{n}{l}^{-1},$

\item  $h(m,\lambda'({l})) = ~({}^{m}l)~{l}^{-1},$

\item  $h(mm_1, n) = ~{}^{m}h(m_1,n)h(m,n),$ 

\item  $h(m,{n}n_1)=h(m,{n})~{}^{{n}}h(m,n_1),$

\item  $h(~{}^{p}m, ~ {}^{p}n) = {}^{p}h(m, n),$
\end{enumerate}
The category of crossed squares will be denoted, 
$\mathfrak{Crs^2}.$\\

In the simplest examples of crossed squares (see \cite{porter}), $\mu$ and
${\mu'}$ are normal subgroup inclusions and $L = M \cap N,$ with $h$ being the conjugation map. We also note that if 
$$
\diagram
{\bf M}\cap{\bf N}\rto\dto&{\bf M}\dto\\
{\bf N}\rto&{\bf G}
\enddiagram
$$
is a simplicial crossed square constructed from a simplicial group ${\bf G}$
and two simplicial normal subgroups ${\bf M}$ and ${\bf N}$ then applying
$\pi_0$, the square gives a crossed square and that up to isomorphism all
crossed squares arise in this way, again see \cite{porter}.

Although when first defined by D. Guin-Walery and J.-L. Loday \cite{wl}, the notion of crossed squares was not linked to that of cat$^2$-groups, it was in this form that Loday gave their generalisation to an $n$-fold structure, cat$^n$-groups (see \cite{loday}).

Recall from \cite{loday} that  a cat$^1$-group is a triple $(G,s,t)$, where
$G$ is a group and $s,t$ are endomorphisms of $G$ satisfying conditions\\
(i)\quad $st=t$ and $ts=s.$\\
(ii)\quad $[\mbox{Ker}s, \mbox{Ker}t] =1.$

It was shown \cite{loday} that setting $M=\mbox{Ker}s,$ $N=\mbox{Im}s$ and
$\partial=t|M,$ then the action of $N$ on $M$ by conjugation  within $G$ makes
$\partial:M\to N$ into a crossed module. Conversely if $\partial: M\to N$ is a
crossed module, then setting $G = M\rtimes N$ and letting $s,t$ be defined by 
$$
s(m,n)=(1,n)
$$
and
$$
t(m,n)=(1,\partial(m)n)
$$
for $m\in M,~n\in N,$ we have that $(G,s,t)$ is a cat$^1$-group.

For a cat$^2$-group, we again have a group, $G$, but this time with two independent cat$^1$-group structures on it. Explicitly:

A cat$^2$-group is a 5-tuple $(G,s_1,t_1,s_2,t_2),$ where $(G,s_i,t_i),$ ~~$i=1,2,$ are cat$^1$-groups and 
$$
s_is_j=s_js_i,\quad t_it_j=t_jt_i,\quad s_it_j=t_js_i
$$ 
for $i,j=1,2$, \quad $i\ne j.$
\begin{thm}{\rm\cite{loday}}
There is an equivalence of categories between the category of cat$^2$-groups and that of crossed squares.
\end{thm}
We include a sketch of the proof as it contains ideas that will be needed
later.

\begin{pf}
The cat$^1$-group $(G,s_1,t_1)$ will gives us a crossed module with
$M=\mbox{Ker}s,~N=\mbox{Im}s$, and $\partial= t|M,$ but as the two cat$^1$-group structures are independent, $(G,s_2,t_2)$ restricts to give cat$^1$-group structures on $M$ and $N$ and makes $\partial$ a morphism of cat$^1$-groups. We thus get a morphism of crossed modules
$$
\diagram
\mbox{Ker}s_1\cap\mbox{Ker}s_2\rto\dto &\mbox{Im}s_1\cap\mbox{Ker}s_2\dto\\
\mbox{Ker}s_2\cap\mbox{Im}s_1\rto&\mbox{Im}s_1\cap\mbox{Im}s_2
\enddiagram
$$
where each morphism is a crossed module for the natural action, i.e.,
conjugation in $G.$ It remains to produce an $h$-map, but this is given by the
commutator within $G$ since if $x\in\mbox{Ker}s_2\cap\mbox{Im}s_1$ and
$y\in\mbox{Im}s_2\cap\mbox{Ker}s_1$ then $[x,y] \in
\mbox{Ker}s_1\cap\mbox{Ker}s_2.$ It is easy to check the axioms for a  crossed square.

Conversely, if 
$$
\diagram
L\rto\dto&M\dto\\
N\rto&P
\enddiagram
$$
is a crossed square, then we can think of it as a morphism of crossed modules
$$
\xymatrixcolsep{0.1pc}
\xymatrixrowsep{0.4pc}
\diagram
L\ddto&&M \ddto\\
& \longrightarrow \\
N&&P.
\enddiagram
$$
Using the equivalence between crossed modules and cat$^1$-groups this gives a morphism
$$
\partial:(L\rtimes N,s,t)\to(M\rtimes P, {s'},{t'})
$$
of cat$^1$-groups. There is an action of $(m,p)\in M\rtimes P$ on $(l,n)\in L\rtimes N$ given by 
$${}^{(m,p)}(l,n) = {}^m({}^pl,{}^pn) = ({}^{\mu(m)p}l
h(m,{}^pn),{}^pn).$$
Using this action, we thus form the associated cat$^1$-group with `big' group $(L\rtimes N)\rtimes (M\rtimes P)$ and induced endomorphisms, $s_1,t_1,s_2,t_2.$
\hfill\end{pf}

\medskip

It is easy to show that cat$^1$-groups are merely a reformulation of an
internal groupoid in the category  $\mathfrak{Grps}$ of groups, whilst
 cat$^2$-groups correspond similarly to double groupoid objects in $\mathfrak{Grps}$.
\section{Crossed $n$-cubes and simplicial groups}
In the form given above, crossed squares were difficult to generalise to higher 
order structures, although cat$^2$-groups could clearly and easily be
generalised to cat$^n$-groups.
The following generalisation is due to Ellis and Steiner \cite{es} and
includes a reformulation of crossed squares as a special case.. Let  $\langle n \rangle$ denote the set $\{1,...,n\}.$

{\rm A crossed $n$-cube of group} is a family $ \{ \mathfrak{M}_A : A\subseteq \langle n \rangle\}$  
of groups, together with homomorphisms $\mu _i:\mathfrak{M}_A\longrightarrow 
\mathfrak{M}_{A\setminus\{i\}}$ for $i\in \langle n \rangle$ 
and functions 
$$h:\mathfrak{M}_A\times \mathfrak{M}_B\longrightarrow \mathfrak{M}_{A\cup B}$$ 
for $A,B\subseteq \langle n \rangle,$ such that if ${}^{a}b$ denotes $h(a,b)b$ for $a \in 
\mathfrak{M}_{A}$ and $b \in \mathfrak{M}_{B}$
with $A\subseteq B,$ then for $a, {a'} \in \mathfrak{M}_{A}$ and $b,{b'} \in 
\mathfrak{M}_{B}, c \in \mathfrak{M}_{C}$ and 
$i,j \in \langle n \rangle,$ the following axioms hold: 
$$
\begin{array}{ll}
1) & \mu _ia = a\ \quad 
\text{{\rm if}}\ i\not \in A, \\ 
2) & \mu _i\mu _ja = \mu _j\mu _ia, \\ 
3) & \mu _ih(a, ~b) = h(\mu _ia, ~\mu _ib), \\ 
4) & h(a, ~b) = h(\mu _ia, ~b)=h(a, ~\mu _ib),  \quad 
\text{{\rm if}}\ i\in A\cap B, \\ 
5) & h(a, ~a^{\prime }) =\lbrack a,{~}a^{\prime }
\rbrack, \\ 
6) & h(a, ~b) = h{(b, ~a)}^{-1}, \\ 
7) & h(a, ~b)=1, \quad\text{if $a = 1$ ~~\text{or}~~ $b = 1,$}\\
8) & h(aa^{\prime }, ~b)={}^{a} h(a^{\prime }, ~b)h(a, ~b), \\ 
9) & h(a,~bb^{\prime })= h(a,~b){~}{}^{b}h(a, ~b^{\prime }), \\ 
10) & {}^{a}h(b, ~c)=
h({}^{a}b,{~}{}^{a}c), \quad\text{{\rm if }} A \subseteq B\cap C,\\ 
11)&{}^{a}h(h(a^{-1}, ~b), ~c)~{}^{c}h(h(c^{-1}, ~a), ~b)~{}^{b}h(h(b^{-1}, ~c), ~a) = 1. \\  
\end{array}
$$

{\em A morphism of crossed n-cubes} is defined in the obvious way: It is a family of group homomorphisms, for $A\subseteq \langle n \rangle,$ 
$
f_A:\mathfrak{M}_A\longrightarrow \mathfrak{M}_{A'}
$
commuting with the $\mu _i$'s and $h$'s. We thus obtain a category of
crossed $n$-cubes denoted by $\mathfrak{Crs^n},$ cf. Ellis and Steiner
\cite{es}.\\
We will concentrate most attention on crossed modules and crossed squares, but 
will recall some of the general theory.

{\noindent\bf Examples:} (1)\quad
For $n=1,$ a crossed $1$-cube is the same as a crossed module.\\
(2)\quad For $n=2,$ one has a crossed square:\\
$$
\diagram
\mathfrak{M}_{\langle 2 \rangle} \dto_{\mu_1} \rto^{\mu_2} & \mathfrak{M}_{\{1\}} \dto^{\mu_1} \\
\mathfrak{M}_{\{2\}} \rto_{\mu_2} & \mathfrak{M}_{\emptyset}.
\enddiagram
$$
Each $\mu _i$ is a crossed module, as is $\mu _1\mu _2$. The h-functions give actions and a function 
$$
h:\mathfrak{M}_{\{1\}}\times \mathfrak{M}_{\{2\}}\longrightarrow \mathfrak{M}_{\langle 2 \rangle}. 
$$
The maps $\mu _2$ (or $\mu _1)$ also define a map of crossed modules. In fact a crossed square can be thought of as a crossed module in the category of crossed modules.\\
(3)\quad
Let {\bf G} be a simplicial group. Then the following
diagram, which will be denoted $\mathfrak{{M}}({\bf G},2),$ 
$$
\diagram
NG_2/\partial_3NG_3 \dto_{~~~~~\partial_2 '} \rto^{\qquad \partial_2} & NG_1\dto^{\mu} \\
\overline{NG_1} \rto_{~~~\mu'} & G_1 
\enddiagram
$$
is a crossed square. Here $NG_1=${\rm Ker}$d_0^1$ and $\overline{NG}_1=${\rm %
Ker}$d_1^1$.

Since $G_1$ acts on $NG_2/\partial _3NG_3,\ \overline{NG}_1$ and $NG_1,$
there are actions of $\overline{NG}_1$ on $NG_2/\partial _3NG_3$ and $NG_1$
via ${\mu'},$ and $NG_1$ acts on $NG_2/\partial _3NG_3$ and $\overline{NG}%
_1$ via $\mu.$ Both $\mu $ and ${\mu'}$ are
inclusions, and all actions are given by conjugation. The h-map is 
$$
\begin{array}{ccl}
NG_1\times \overline{NG}_1 & \longrightarrow & NG_2/\partial _3NG_3 \\ 
(x,\overline{y}) & \longmapsto  & h(x,~y)= \lbrack s_1x,~s_1ys_0{y}^{-1} \rbrack \partial_3NG_3.
\end{array}
$$
Here $x$ and $y$ are both in $NG_1$ as there exists a
bijection between $NG_1$ and $\overline{NG}_1$,  the element $\overline{y}$ is
the image of $y$ under this. 

This last example effectively presents a functor
$$
\diagram
\mathfrak{{M}}~:~\mathfrak{SimpGrp} \rto & \mathfrak{Crs^2}.
\enddiagram
$$ 
(4)\quad Let $G$ be a group with normal subgroups ${N}_1,\ldots ,{N}_n$. Let 
$$
\begin{array}{ccc}
\mathfrak{M}_A=\bigcap \{{N}_i:i\in A\} & \text{and} & \mathfrak{M}_\emptyset =${G}$
\end{array}
$$
with $A\subseteq \langle n \rangle.$ For $i\in \langle n \rangle,$   $\mathfrak{M}_A$ is a normal  
subgroup of $\mathfrak{M}_{A\setminus\{i\}.}$ Define 
$$
\mu _i:\mathfrak{M}_A\longrightarrow \mathfrak{M}_{A\setminus\{i\}} 
$$
to be the inclusion. If $A,B\subseteq \langle n \rangle$, then $\mathfrak{M}_{A\cup B}=\mathfrak{M}_A\cap \mathfrak{M}_B,$
let 
$$
\begin{array}{cccl}
h: & \mathfrak{M}_A\times\mathfrak{M}_B & \longrightarrow  & \mathfrak{M}_{A\cup B} \\  
& (a,b) & \longmapsto  & [{~}a,{~}b{~}]
\end{array}
$$
as $[\mathfrak{M}_A, \mathfrak{M}_B] \subseteq \mathfrak{M}_A\cap \mathfrak{M}_B,$ 
where $a\in \mathfrak{M}_A,\ b\in \mathfrak{M}_B.$ Then 
$$
\{\mathfrak{M}_A:\ A\subseteq \langle n \rangle,\ \mu _i,\ h\} 
$$
is a crossed $n$-cube, called the {\em inclusion crossed n-cube} given by the
normal $n$-ad of groups $(G;\ {N}_1,\ldots ,{N}_n).$  The following result is
then fairly easily proved, see \cite{porter}.

\begin{prop}\label{y}
Let $({\bf G};\ {N}_1,\ldots ,{N}_n)$ be a simplicial normal $n$-ad of 
subgroups of groups and define for $A\subseteq \langle n \rangle$%
$$
\mathfrak{M}_A=\pi _0(\bigcap\limits_{i\in A}{N}_i) 
$$
with homomorphisms $\mu _i:\mathfrak{M}_A\longrightarrow \mathfrak{M}_{A\setminus\{i\}}$ and h-maps induced by
the corresponding maps in the simplicial inclusion crossed $n$-cube,
constructed by applying the previous example to each level. Then $\{\mathfrak{M}_A:\
A\subseteq \langle n \rangle,\ \mu _i,\ h\}$ is a crossed $n$-cube.\hfill $\square$
\end{prop}

 Up to isomorphism, all crossed $n$-cubes arise in this way.
 In fact any crossed $n$-cube can be realised (up to isomorphism) as a $\pi _0$
 of a simplicial inclusion crossed $n$-cube coming from a simplicial normal 
 $n$-ad of groups. 
 
In 1993, the second author, \cite{porter}, described  a functor from the
category of simplicial groups to that of crossed $n$-cubes of groups. We will summarise
its construction.

The functor  is constructed using the \emph{d\'ecalage} functor studied
by Duskin \cite{duskin} and Illusie \cite{illusie} and is a $\pi _0$-image of 
a functor taking values in a category of simplicial normal $(n+1)$-ads. 
The \emph{d\'ecalage} functor will be denoted
by $\mathfrak{Dec}.$ Given any simplicial group ${\bf G},$ 
$\mathfrak{Dec}{\bf G}$ is the augmented simplicial group obtained from 
${\bf G}$ by forgetting the zeroth face and degeneracy operators at each level
and then renumbering the levels (cf. Duskin \cite{duskin} or Illusie
\cite{illusie}). (There are two main forms of the \emph{d\'ecalage} functor.  The
alternative forgets the last face instead of the zeroth one.  This second form 
is used in \cite{duskin}.  Reversing the indexed order of the faces used
gives an equivalence between the two theories and allows for the easy translation of proofs between them.)
In the convention adopted in this paper, we thus have $$\mathfrak{Dec}{ G}_n = G_{n+1},$$ with face operators
$$d_i^{n,Dec} = d_{i+1}^{n+1} $$ and degenerate operators 
$$s_i^{n,Dec} = s_{i+1}^{n+1}.$$
The remaining degeneracy $s^{n+1}_0$ of ${\bf G}$  yields 
a contraction of $\mathfrak{Dec}^1{\bf G}$ as an augmented simplicial group and
$$
{\rm \mathfrak{Dec}}^1{\bf G}\simeq {\bf K}(G_0,0), 
$$  
by an explicit natural homotopy equivalence (cf. Duskin \cite{duskin}). The 
zeroth face map will be denoted $\delta_0 :{\mathfrak{Dec}}^1{\bf G}
\longrightarrow {\bf G}$.  This is a split epimorphism and has kernel the 
simplicial group, Ker~$d_{0}$ used above. We have
$\pi_0(\mbox{Ker}\delta) \longrightarrow \pi_0(\mathfrak{Dec}{\bf G})$ is a
crossed module and using that
$$\pi_0 {\bf H} = NH_0/\partial NH_1$$
for any simplicial group ${\bf H}$, we get: \\
{\bf Case 1:}\qquad ${\bf H} = \mathfrak{Dec}{\bf G}.$ \\
Then $NH_0 = G_1,$  $NH_1 = $Ker$(d_1^2: G_2 \longrightarrow G_1)$ so
$$
\begin{array}{rcl}
\pi_0(\mathfrak{Dec}{\bf G}) & \cong & G_1/d_2^2(\text{Ker}d_1^2), \\
& \cong & G_0,
\end{array}
$$
since $G_1 \cong$ Ker$d_1^1\rtimes s_0(G_0).$ \\
{\bf Case 2:}\qquad ${\bf H} = $ Ker$\delta.$  Then $NH_0 = $ Ker$d_0^1 = NG_1,$ 
$NH_1 = $  Ker$d_0^2 ~\cap $ Ker$d_1^2 = NG_2$ 
so $$\pi_0(\textrm{ Ker }\delta) = NG_1/\partial_2NG_2.$$

\medskip

Iterating the $\mathfrak{Dec}$ construction gives an augmented bisimplicial 
group%
$$
\diagram
(...~{\rm \mathfrak{Dec}}^3{\bf G}\rto<.9ex> \rto \rto<-.9ex> &{\rm \mathfrak{Dec}}^2{\bf G}
\rto<.5ex>^{\textstyle \delta_0 } \rto<-.5ex>_{\textstyle \delta_1}
 &{\rm \mathfrak{Dec}}^1{\bf G}) 
\enddiagram
$$
which in expanded form is the total \emph{d\'ecalage} of {\bf G},
(see \cite{duskin} or \cite{illusie} for details). The maps from $\mathfrak{Dec}^i{\bf G}$
to $\mathfrak{Dec}^{i-1}{\bf G}$ coming from the $i$ first face maps will be labelled $
\delta_0,\ldots ,\delta_{i-1}$ so that
$\delta _0= d_{0},\,\delta_1= d_{1}$ and so on.

For a simplicial group {\bf G} and a given $n,$ we write 
$\mathfrak{M}({\bf G},n)$ for the crossed $n$-cube arising as a functor 
$$
{\mathfrak{M}(-,n):\mathfrak{SimpGrp} \longrightarrow \mathfrak{Crs^n.}} 
$$
which is given by 
$\pi_0(\mathfrak{Dec}{\bf G}; \mbox{Ker}\delta_0,\ldots,\mbox{Ker}\delta_{n-1}).$
The following data explicitly gives this crossed $n$-cube of groups,  for the details see \cite{porter}:
\begin{thm}\label{t1}
If {\bf G} is a simplicial group, then the crossed $n$-cube $\mathfrak{M}${\rm (}%
{\bf G},$n${\rm )} is determined by:

(i) for $A\subseteq \langle n \rangle,$%
$$
\mathfrak{M}({\bf G},n)_A=\frac{\bigcap_{j\in A}\text{{\rm Ker}}d_{j-1}^n}{%
d_{n+1}^{n+1}(\text{{\rm Ker}}d_0^{n+1}\cap \{\bigcap_{j\in A}\text{{\rm Ker}%
}d_j^{n+1}\})}; 
$$

(ii) the inclusion 
$$
\bigcap_{j\in A}\text{{\rm Ker}}d_{j-1}^n\longrightarrow \bigcap_{j\in
A\setminus\{i\}}\text{{\rm Ker}}d_{j-1}^n 
$$
induces the morphism 
$$
\mu _i:\mathfrak{M}({\bf G},n)_A\longrightarrow \mathfrak{M}({\bf G},n)_{A\setminus\{i\}}; 
$$

(iii) the functions, for $A,B\subseteq \langle n \rangle,$ 
$$
h:\mathfrak{M}({\bf G},n)_A\times \mathfrak{M}({\bf G},n)_B\longrightarrow 
\mathfrak{M}({\bf G},n)_{A\cup B} 
$$
are given by 
$$
h(\bar x,\bar y) = \overline{[x,~y]}, 
$$
where an element of $\mathfrak{M}({\bf G},n)_A$ is denoted by $\bar{x}$ with $x\in
\bigcap_{j\in A}${\rm Ker}$d_{j-1}^n.$\hfill $\square$
\end{thm}
 In general, we use the
($n-1$)-skeleton of the total \emph{d\'ecalage} to form an $n$-cube and thus a 
simplicial inclusion crossed $n$-cube of kernels. Continuing this $n$-times gives 
the simplicial inclusion crossed $n$-cube corresponding to the simplicial normal 
$(n+1)$-ad, 
$$\mathcal{M}({\bf G},n)=(\mathfrak{Dec}^n{\bf G};{~}\text{Ker}\delta_{0},\ldots ,\text{Ker}\delta_{n-1}),$$
and its associated crossed $n$-cube is
$$
\pi _0(\mathcal{M}({\bf G},n))=\mathfrak{M}({\bf G},n). 
$$
The results of \cite{porter} now follow by direct calculation on examining the construction of
$\pi _0$ as the zeroth homology of the Moore complex of each term in the inclusion
crossed $n$-cube, $\mathcal{M}({\bf G}, n).$

Expanding this out for low values of $n$ gives:\\
\noindent 1) For $n=0$, 
$$
\begin{array}{rcl}
\mathfrak{M}({\bf G},0) & = & G_0/d_1(
\text{{\rm Ker}}d_0,) \\  & \cong  & \pi _0(
{\bf G}), \\  & = & H_0({\bf G).}
\end{array}
$$
2) For $n=1$, $\mathfrak{M}({\bf G},1)$ is the crossed module 
$$
\mu_1 :\text{{\rm Ker}}d_0^1/d_2^2(NG_2)\longrightarrow G_1/d_2^2(\text{{\rm Ker}}d_0^2). 
$$
Since $d_2^2(NG_2)=$ [{\rm Ker}$d_1^1,${~}{\rm Ker}$d_0^1],$ this gives 
$$
\mu :NG_1/[\text{{\rm Ker}}d_1^1,{~}\text{{\rm Ker}}d_0^1]\longrightarrow
G_0. $$
$$
\mathfrak{M}(\mathbf{G}, 1) \cong
\left( \begin{array}{ccc}
NG_1/\partial_2NG_2 \longrightarrow G_0
\end{array}
\right).
$$

3) For $n=2$, $\mathfrak{M}({\bf G},2)$ is 
$$
\diagram
\mbox{\rm Ker}d^{2}_{0} \cap \mbox{\rm Ker}d^{2}_{1} / 
d^{3}_{3} (\mbox{ \rm Ker}d^{3}_{0} 
\cap \mbox{\rm Ker}d^{3}_{1} \cap \mbox{\rm Ker}d^{3}_{2})\dto_{\quad\mu_1} 
\rto^{\qquad\quad \mu_2} & \mbox{\rm Ker}d^{2}_{0}  /  d^{3}_{3}(\mbox{\rm Ker}d^{3}_{0} 
\cap \mbox{\rm Ker}d^{3}_{1}) \dto^{\mu_1} \\
\mbox{\rm Ker}d^{2}_{1}  /  d^{3}_{3}(\mbox{\rm  Ker }d^{3}_{0} 
\cap \mbox{\rm Ker}d^{3}_{2}) \rto_{~~~~~\mu_2}  & G_2 / d^{3}_{3}
(\mbox{\rm Ker}d^{3}_{0}).
\enddiagram
$$
As shown in \cite{porter}, this is isomorphic to 
$$
\diagram
NG_2 / d^{3}_{3} (NG_3 ) \dto_{\mu_1} \rto^{\quad\mu_2}   & \mbox{\rm Ker}d^{1}_{0} \dto^{\mu_1} \\
\mbox{\rm Ker}d^{1}_{1} \rto_{\mu_2} &  G_1,   
\enddiagram
$$
that is
$$
\mathfrak{M}(\mathbf{G}, 2) \cong
\left ( \diagram
NG_2 / \partial_3 (NG_3 ) \dto \rto & \mbox{ Ker }d_{0} 
\dto \\
\mbox{ Ker }d_{1} \rto  &   G_1    
\enddiagram \right).
$$
 Here the $h$-map is 
$$h: \mbox{Ker}d_0^1\times\mbox{Ker}d_1^1\longrightarrow NG_2/d_3^3(NG_3)$$
given by $h(x,y) = [s_1x, ~s_1ys_0y^{-1}]~\partial_3NG_3.$
Note if we consider the above crossed square as a vertical morphism
of crossed modules we can take its kernel and cokernel within the category of crossed
modules. In the above, the morphisms in the top left hand corner are 
induced from $d_2$ so
$$
\mbox{Ker}\left ( 
\mu_1 : \frac{NG_2}{\partial_3NG_3}\longrightarrow  \mbox{Ker}d_1
\right) =  \frac{NG_2\cap \mbox{Ker}d_2}{\partial_3NG_3} 
 \cong  \pi_2({\bf G})
$$
whilst the other map labelled $\mu_1$ is an inclusion so has trivial kernel.
Hence the kernel of this morphism of crossed modules is
$$
\pi_2({\bf G})\longrightarrow 1.
$$
The image of $\mu_2$ is  normal in both the simplicial groups 
on the bottom line and as $\mbox{Ker}d_0 =NG_1$ with the corresponding $\mbox{Im}\mu_1$
being $d_2NG_2,$ the cokernel is $NG_1/\partial_2NG_2,$ whilst 
$G_1/\mbox{Ker}d_0\cong G_0,$ i.e., the cokernel of $\mu_1$ is 
$\mathfrak{M}({\bf G}, 1).$

In fact of course $\mu_1$ is not only a morphism of crossed modules, it is 
a crossed module. This means that $\pi_2({\bf G})\longrightarrow 1$ is in 
some sense a ~$\mathfrak{M}({\bf G}, 1)$-module and that~ 
$\mathfrak{M}({\bf G}, 2)$ can be thought of as a crossed extension of~ 
$\mathfrak{M}({\bf G}, 1)$ by $\pi_2({\bf G}).$
\section{2-crossed modules from crossed squares}
D. Conduch\'e's unpublished work shows that there exists an equivalence (up to homotopy) between the category of crossed squares of groups and that
of 2-crossed modules of groups. Loday defined a mapping complex of crossed squares by: if
$${\cal M} = 
\diagram
{L}\dto_{\lambda'} \rto^{\lambda} & M\dto^{\mu} \\
{N}  \rto_{\nu} & {P} 
\enddiagram
$$
is a crossed square, then its mapping complex is
$$
\diagram
{L} \rrto^{\partial_2} && M\rtimes{N}  \rrto^{\partial_1}&&{P}
\enddiagram
$$
where $\partial_2l = (\lambda{l}^{-1}, ~{\lambda'}l)$ and 
$\partial_1(m, ~{n}) = \mu(m){\nu}({n}).$ 

Conduch\'{e} showed that this mapping complex is always a 2-crossed complex,
representing the same 3-type as the original crossed square. Why?

\subsection{From crossed squares to bisimplicial groups}
Recall that from a crossed module ${\cal M} = (M \stackrel{\mu}{\to} P)$, we can build a
simplicial group whose Moore complex is trivial in dimensions 2 and above.
This is the nerve of the associated cat$^1$-group. (Because a cat$^1$-group is 
an internal groupoid in $\mathfrak{Grps}$, we can form the nerve of its
category structure internally within $\mathfrak{Grps}$ and hence obtain a
simplicial group.)  We need an explicit description of this, elementwise:  The 
simplicial group $Ner({\cal M})$ has 
$$Ner({\cal M})_0 = P$$
$$Ner({\cal M})_n = M\rtimes ( \ldots (M\rtimes P) \ldots),$$
with $n$ semidirect factors of $M$ (see Conduch\'{e}, \cite{conduche},
Carrasco and Cegarra, \cite{C&C1}, or our own papers \cite{mp1,mp5} for more
on  semidirect decompositions and simplicial groups).  If $(m,p) \in
Ner({\cal M})_1$, then 
$$d_0^1(m,p) = \mu(m)p$$
$$d^1_1(m,p) = p$$
and $$s^0_0(p) = (1,p).$$
If $(m_2,m_1,p) \in
Ner({\cal M})_2$, then 
$$d_0^2(m_2,m_1,p) = (m_2,\mu(m_1)p)$$
$$d_1^2(m_2,m_1,p) = (m_2m_1,p)$$
$$d_1^2(m_2,m_1,p) = (m_1,p)$$
$$s^1_0(m,p) = (m,1,p)$$
$$s^1_1(m,p) = (1,m,p)$$
and the obvious pattern continues to higher dimensions.

Now going to crossed squares, let 
$${\cal M} = 
\diagram
{L}\dto_{\lambda'} \rto^{\lambda} & M\dto^{\mu} \\
{N}  \rto_{\mu'} & {P} 
\enddiagram
$$
be a crossed square.  The construction we will need is the bisimplicial
nerve of ${\cal M}$ or rather of the associated cat$^2$-group.  That will be
an internal double groupoid in $\mathfrak{Grps}$ and so when we take the
nerves in the two directions, we will get a bisimplicial group.  Although
again this is known, we will want explicit descriptions of elements etc. for
explicit calculations of actions, pairings, etc. later on, so will give quite
a lot of detail, even when it is fairly simple to check.

Considering the crossed square ${\cal M}$ as a morphism /crossed module from
$(L\to N)$ to $(M\to P)$, we apply the above nerve construction to its domain
and codomain to get a crossed module of simplicial groups.  In low dimensions:
$$
\xymatrix{L\rtimes L\rtimes N \ar[r]^{(\lambda,\lambda,\nu)}
  \ar[d]<-1ex>\ar[d]\ar[d]<1ex> & M\rtimes M\rtimes P\ar[d]<-1ex>\ar[d]\ar[d]<1ex> \\
  L\rtimes N \ar[r]^{(\lambda,\nu)}\ar[d]<-.5ex>\ar[d]<.5ex>& M\rtimes P\ar[d]<-.5ex>\ar[d]<.5ex>\\
N \ar[r]^{\nu}& P}
$$
To check, for instance, that $(\lambda,\nu)$ is a crossed module, you need to
use the $h$-map of $\cal M$, as follows:  First the action of $(m,p)$ on
$(l,n)$ is $${}^{(m,p)}(l,n) = {}^m({}^pl,{}^pn) = ({}^{\mu(m)p}l
h(m,{}^pn),{}^pn).$$
Now writing $\partial = (\lambda,\nu)$,
$$\partial
({}^{(m,p)}(l,n))=(\lambda({}^{\mu(m)p}l).\lambda(h(m,{}^pn)),\nu({}^pn),$$but 
$\lambda(h(m,{}^pn))= m.{}^{\nu{{}^pn}}m^{-1}$ and so a routine calculation
shows that this expands to $$(m,p)(\lambda l,\nu n)({}^{p^{-1}}m^{-1}, p^{-1}),$$ i.e. to 
${}^{(m,p)}\lambda(l,n)$, so $\partial$ satisfies the first crossed module
axiom.  We leave the second crossed module axiom (Peiffer identity) to the
diligent reader.  It uses $h(\lambda l, n^\prime) = l {}^{n^\prime}l^{-1}$.
The description of the `vertical' face and degeneracies is as before.  

Next we 
start building the nerve in the second direction. Writing $\cal X$ for
the resulting simplicial group, we get:
\begin{eqnarray*}
{\cal X}_{0,0} & = & P\\
{\cal X}_{0,1} & = & M\rtimes P\\
\mbox{and in general\rule{5.45cm}{0cm} }  & & \rule{5cm}{0cm}\\
{\cal X}_{0,q} & = &M\rtimes\ldots M\rtimes P = M^{(q)}\rtimes P\rule{5cm}{0cm}\\
\mbox{with $q$-factors of $M$, (note \rule{3.2cm}{0cm} }&& \\
\mbox{the shorthand version).\rule{3.95cm}{0cm}} &&\\
\mbox{Similarly\rule{5.7cm}{0cm} }{\cal X}_{p,0} & = &N\rtimes\ldots N\rtimes
P =
N^{(p)}\rtimes P
\end{eqnarray*}
with $p$-factors of $N$. 
In general ${\cal X}_{p,q}$ can be written :
$${\cal X}_{p,q} = (L^{(q)}\rtimes N)^{(p)}\rtimes(M^{(q)}\rtimes P).$$

\textbf{Remark:}

This initially looks asymmetric but in fact is not.  The result of forming the 
horizontal nerve first then the vertical one would give an isomorphic
bisimplicial group.  This is easy to show categorically, but is tedious to
show elementwise as it makes repeated use of the interchange /$h$-map
structure within ${\cal M}$.

\subsection{From bisimplicial groups to simplicial groups.}
  There are two useful ways of passing from  bisimplicial groups to simplicial
  groups.  One is the diagonal, the other, due to Artin and Mazur, \cite{am},
  is the `codiagonal' and, for us, is more useful.  In the bisimplicial group ${\cal X}$, the corresponding 
  Moore bicomplex has relatively few non-zero terms.  In fact, $N({\cal
    X})_{p,q} $ will be zero if $p$ or $q$ is bigger than 1.  If you take the
  diagonal and try to apply the Moore complex functor, the task looks
  horrendous. Although, in fact, this Moore complex has length 2, the individual 
  terms are quite large with a complicated expression for the differential.
  The diagonal gets complicated because its $p$-simplexes are the
  $(p,p)$-simplices of the original bisimplicial group and so seem to
  correspond to $2p$-dimensional data. (In fact they need a list of $(p+1)^2$-elements to describe them.)

In the Artin-Mazur construction, if ${\cal X}_{*,*}$ is a bisimplicial group,
we first form for each $n$, a group
$${\cal X}_{(n)} = \prod_{p + q = n} {\cal X}_{p,q}.$$
Within this ${\cal X}_{(n)}$, we pick out a subgroup, $\nabla({\cal X})_{n}$
as follows:
Let $\underline{x} = (x_0, \ldots, x_n) \in {\cal X}_{(n)}$ with $x_p \in {\cal
  X}_{p,n-p}$.  Then $\underline{x} \in \nabla({\cal X})_{n}$ if and only if for each
$p = 0, \ldots, n-1$
$$d_0^vx_p = d^h_{p+1}x_{p+1}.$$
This mysterious formula can be best remembered by looking at the case $n = 2$
and the diagram
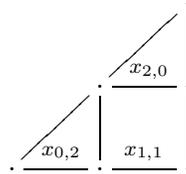
\vspace{-2cm}
\begin{center}
\rule{0mm}{2cm}\xymatrix{ & & .\ar@{-}[d]\\
 & . \ar@{-}[ur]\ar@{-}[r]^{~x_{2,0}}&.\\
 . \ar@{-}[ur]\ar@{-}[r]^{~x_{0,2}}&.\ar@{-}[r]^{x_{1,1}}\ar@{-}[u]&.\ar@{-}[u]}\begin{eqnarray*}d_0^vx_0 = d^h_{1}x_{1}\\
d_0^vx_1 = d^h_{2}x_{2}
\end{eqnarray*}
\end{center}
where we have expanded the notation $x_p$ to $x_{p,n-p}$ to make the changes
in dimension clearer (we hope).  Similar subdivided simplex diagrams can be
seen to give the formulae in higher dimensions, although above $n=3$, they
cannot be so simply drawn.  The link between this and the ordinal subdivision
of the paper, \cite{ep2}, by Ehlers and the second author will be explored and 
exploited later.

The face and degeneracy maps of $\nabla({\cal X})$ are built up in an obvious
way from this formula:
$$d_j = d_j^\nabla :\nabla({\cal X})_n\to \nabla({\cal X})_{n-1},$$
for $\underline{x} = (x_0, \ldots, x_n),$ with $d_0^vx_p = d_{p+1}^hx_{p+1},$
$p = 0, \ldots, n-1$, then for $0 < j < n$,\\
$d_j^\nabla(\underline{x}) := (d_j^vx_0,d_{j-1}^vx_1,
\ldots,  d_1^vx_{j-1}, d_j^hx_{j+1}, \ldots,  d_j^hx_n)$;\\
$d_0^\nabla(\underline{x}) := (d_0^hx_1,d_0^hx_2, \ldots,d_0^hx_n)$;\\
and\\
$d_n^\nabla(\underline{x}) := (d_n^vx_0,d_{n-1}^vx_1,
\ldots, d_1^vx_{n-1}),$\\
whilst
$s_i^\nabla(\underline{x}) := (s_i^vx_0,s_{i-1}^vx_1, \ldots , s_0^vx_i,
s_i^hx_i, \ldots, s_i^hx_n)$ for $0 \leq i \leq n$.

\medskip

\textbf{Remark:}

Of course the formulae for $d_0$ and $d_n$ can be considered as being special
cases of that for $d_j$.

\medskip

Suppose now given a bisimplicial group, ${\cal X}$, such  that for any fixed
$p$, the Moore complex of ${\cal X}_{p, *}$ is of length at most  1 (so is a
crossed module) and that for any  fixed $q$, similarly $N({\cal X}_{*,q})$ has 
length at most 1.  Then we claim that $N(\nabla({\cal X}))$ has length $\leq
2$.

We examine an element $\underline{x} \in N(\nabla({\cal X}))_n$ for $n \geq
3$.  Since $d_0(\underline{x})= 1$, we have 
$$x_p \in Ker d^h_0 \quad \mbox{ for all } p > 0.$$
Now for each $0 < j < n$, $d_j(\underline{x})= 1$, so
$$d_j^vx_0 = d_{j-1}^vx_1 = \ldots d_1^vx_{j-1} = 1,$$
whilst 
$$d_j^hx_{j+1} = \ldots = d_j^hx_n = 1.$$
Looking at $x_n$, $d_j^hx_n = 1$ for $j = 0, \ldots, n-1$, i.e. $x_n \in
N({\cal X}_{*,0})_n$.  If $n \geq 2$, then this group is trivial, so $x_n =1$, 
but $d^h_nx_n = d^v_0x_{n-1}$, so $d^v_0x_{n-1}=1$ as well.  Turning to
$x_{n-1}\in {\cal X}_{n-1,1}$, we have $d_0^hx_{n-1}=1$. From $d_j^\nabla$, we 
have $d_j^hx_{n-1}=1$ for $j = 1, \ldots, n-2$, i.e. $x_{n-1}\in N({\cal
  X}_{*,1})_{n-1}$.  As $n \geq 3$, $n-1\geq2$, so $ N({\cal X}_{*,1})_{n-1} =
1$, i.e. $x_{n-1}=1$ and hence $d_0^v(x_{n-1})=1$.  Continuing like this, we
get $x_k = 1$ for all $k \geq 2$: for each index, we already know 
$$d_j^hx_k=1 \mbox { for } j \leq k-1,$$
i.e. $x_k \in N({\cal X}_{*,n-k})_k = 1 $ if $k \geq 2$.  This leaves us with
just $x_0$ and $x_1$ to examine:
$$d_0^vx_1 = d_2^hx_2 =1,$$
whilst $d_j^vx_1 =1 $ (from $d^\nabla_{j+1}\underline{x} = 1$ so again $x_1
\in N({\cal X}_{1,*})_{n-1}$ which is trivial as $n-1 \geq 2$, i.e. $x_1 = 1$,
so $d_0^vx_0 = 1$  as well.  Finally $d_j^vx_0 = 1$ for all $0 < j < n$ (from
$d^\nabla_{0}\underline{x} = 1$), so $x_0 \in  N({\cal X}_{0,*})_n = 1$.

We have thus proved :
\begin{prop}
If ${\cal X}_{*,*}$ is a bisimplicial group such that for any $p$,  $N({\cal
  X}_{p,*})_q = 1$ for $q \geq 2$ and for any $q$, $N({\cal  X}_{*,q})_p = 1$
for $p \geq 2$, then $N(\nabla ({\cal X})_n=1$ for $n \geq 3$.\hfill$\square$
 \end{prop}

We have in fact proved a bit more.  If ${\cal X}_{*,*}$ satisfies the
conditions of the proposition and $\underline{x} = (x_0,x_1,x_2) \in
N(\nabla({\cal X}))_2$, then $x_2 =1$ since $x_2 \in N({\cal X}_{*,0})_2 =1$.
Of course $N(\nabla({\cal X}))_n = 1$ for $n \geq 3$ implies that
$N(\nabla({\cal X}))$ is a 2-crossed module.  We will examine this in some
detail.

\medskip

\textbf{Remark:}

The above proposition suggests two questions.  Firstly if we weaken the
condition that each direction gives a crossed module as it Moore complex to
that it gives a crossed \emph{complex}, is it true that $N(\nabla({\cal X}))$
is a 2-crossed complex?  Also, is there a generalisation to higher order
multisimplicial groups? We will later on show that the answer to the second
one is positive.

\medskip

For the moment we will look in more detail at the case when ${\cal X}$ is the
binerve of a crossed square.

\subsection{The 2-crossed module structure of $N(\nabla({\cal X}({\cal M})))$.}

As before, we let ${\cal M}$ be a crossed square and will write ${\cal
  X}({\cal M})$ for $Ner^h(Ner^v({\cal M}))$, i.e. for the binerve of ${\cal
  M}$. For convenience we will often write $G$ instead of $\nabla({\cal
  X}({\cal M}))$.

We know $N(G)_n = 1$ if $n \geq 3$.  It is also clear that $G_0 = P$, so
$N(G)_0 = P$ as well.  We thus only need to work out  $N(G)_k$ for $k= 1$ and
$2$, together with explicit formulae for the boundary maps and the Peiffer
pairing from $N(G)_1\times N(G)_1$ to $N(G)_2$.

Suppose $x\in G_1$, then $\underline{x} = (x_0,x_1)$ with $x_0 \in {\cal
  M}_{0,1}$, $x_1 \in {\cal M}_{1,0}$.  We have 
$${\cal M}_{0,1} = M\rtimes P, \hspace{2cm} {\cal M}_{1,0} = N\rtimes P,$$
so $x_0 = (m,p)$, $x_1 = (n,p^\prime)$ and, since $\underline{x} \in
\nabla{\cal X}({\cal M}))_1$
$$d_0^v(m,p) = \mu(m)p = p^\prime = d_1^h(n,p^\prime),$$i.e. $p^\prime$ is
determined by $(m,p)$.  The assignment 
$$\underline{x} = ((m,p), (n,\mu(m)p)) \to (n,m,p)$$ is easily checked to give
an isomorphism
$$G_1\stackrel{\cong}{\to} N\rtimes(M\rtimes P),$$
where $M$ acts on $N$ via $P$, ${}^mn = {}^{\mu(m)}n.$ Identifying $G_1$ with  
$N\rtimes(M\rtimes P)$, $d_0$ and $d_1$ have the descriptions:
\begin{eqnarray*}
d_0(n,m,p) &=& \nu(n)\mu(m)p\\
d_1(n,m,p) &=&p.
\end{eqnarray*}
Thus $\underline{x} \in N(G_1)$  if and only if $p = \mu(m)^{-1}\nu(n)^{-1}$
and it is again easily verified that $$
NG_1 \cong M\rtimes N,$$where the isomorphism is given by
$$(n^{-1},m^{-1},\mu(m)\nu(n)) \to (n,m).$$
We note that via this isomorphism, we get $$\partial_1 : NG_1 \to NG_0$$ is
given by
$$\partial_1(n,m) = \mu(m)\nu(n).$$
  This does look
strange (inverses do not work well with conditions for homomorphisms), but
does work as is easily checked.

Turning to $NG_2$, $\underline{x} \in \nabla{\cal X}_2$ will have the form
$(x_0,x_1,x_2)$ with 
\begin{eqnarray*}
x_0 & =& (m_2,m_1,p) \in M\rtimes(M\rtimes P)\\
x_1 &=& (l,n,m,p^\prime) \in (L\rtimes N)\rtimes (M\rtimes P)\\
x_2 & = & (n_2,n_1,p^{\prime\prime}) \in N\rtimes (N\rtimes P)
\end{eqnarray*}
The equations $d_0^vx_0 = d_1^hx_1$ and $d_0^vx_1 = d_2^hx_2$ give relations
between the individual coordinates implying that $m_2 =m$, $p^\prime =
\mu(m_1)p$, $ \lambda^\prime(l) n = n_1$ and $\mu(m)p^\prime =
p^{\prime\prime}$.  (This is best seen on a diagram as above, but such a
diagram is best left to the reader to draw!) Now suppose $\underline{x} \in
N(\nabla {\cal X})_2$, then in addition one gets $p = 1$, $m_2m_1 = 1$,
etc. and it follows that all the `coordinates' depend only on $l$. In fact, 
$$\underline{x} = ((\lambda (l) ^{-1},(\lambda (l),1)), ((l,\lambda^\prime (l)
^{-1}), (\lambda (l) ^{-1}, \mu \lambda (l))), (1,(1,1)))$$ with
$d_2(\underline{x}) = ((\lambda (l), 1), (\lambda^\prime (l),\mu \lambda (l)))$, 
which identifies to $(\lambda (l)^{-1}, \lambda^\prime (l))\in M\rtimes
N$. Thus the normal chain complex $NG$ is isomorphic to\\
$$\xymatrix{L\ar[r]^{\partial_2\hspace*{3mm}} & M\rtimes N \ar[r]^{\hspace*{4mm}\partial_1} &P}$$ with  
\begin{eqnarray*}
\partial_2(l)& = &((\lambda^\prime (l)^{-1}, \lambda (l)))\\
\partial_1(l)& = & \mu(m)^{-1}\nu(n)^{-1}.
\end{eqnarray*}

To complete the description, we should really specify the Peiffer pairing
$$\{ \quad,\quad \} : NG_1 \times NG_1 \to NG_2.$$
We can use the formula given in our earlier paper, \cite{mp2}, derived from
work of Conduch\'e:
$$\{\underline{x}, \underline{y}\} = s_0(\underline{x})s_1(\underline{y}) s_0(\underline{x})^{-1}s_1(\underline{x}\underline{y}^{-1}\underline{x}^{-1})$$
Since $s_0(\underline{x}) = s_0(x_0,x_1) = (s_0^vx_0,s_0^hx_0,s_0^hx_1)$ and 
$s_1(\underline{x}) = s_1(x_0,x_1) = (s_1^vx_0,s_0^vx_1,s_1^hx_1)$, we can get 
an explicit desciption of $\{\underline{x}, \underline{y}\}$.  In fact, that
description is not that useful.  Its exact form is dependent on the order of
multiplication used.  These forms are equivalent via the relations for
expanding $h(m.m^\prime,n)$ in terms of $h(m,n)$ and $h(m^\prime,n)$, but at
the risk of lengthening the expression. These different forms are fairly
complicated and so have been omitted here, but can be retrieved from the
original formula.

\medskip

\textbf{Remarks}\\
(i) In Conduch\'e's original letter, he gives a much simpler form, namely $h(m, nbn^{-1})$.  Our efforts to reduce the above to
something as simple as this have so far failed!  This is almost certainly due
to the change in convention on the Moore complex, but we have not
managed fully to understand the reason for the greater
complication. Conduch\'e's  elegant form  \emph{does} work as is very easily checked.  What is not
clear is its relationship with the `canonical' form coming from the Peiffer pairings.

(ii) The above shows a subtle difficulty encountered when working with elements in models for homotopy $n$-types for $n > 1$.  With (free) groups, we are used to handling normal forms of elements, but even with free crossed modules, the Peiffer identity makes calculating with representatives of elements much more difficult.  This emphasises the need for a higher order version of both combinatorial and computational group theory, adapting the methods of the classical case to these higher dimensional situations. Some success has been achieved in this area by Alp and Wensley, \cite{MA&CDW}, and by Ellis, \cite{gje2}.

\subsection{Squared complexes and 2-crossed complexes}

In 1993, Ellis defined the notion of a squared complex, \cite{gje}. A crossed
complex combines a crossed module at its `base' with a continuation by a chain 
complex of modules further up.  They thus have good descriptive power, being
able to model 2-types via their crossed module part and also a certain amount
of higher homotopy information, typically thought of as generalising chains on 
the universal cover of a space.

Squared complexes are one of the possible notions that generalise crossed
complexes to include the so-called quadratic information available in a
3-type.  Other versions include double crossed complexes (cf. Tonks,
\cite{andy}), 2-crossed complexes (cf. the authors, \cite{mp4}) and quadratic
complexes, (cf. Baues, \cite{baues}).  We will need to examine 2-crossed complexes in
detail later.
 
A squared complex consists of a diagram of group homomorphisms
$$
\diagram
& & & & N\drto^{\mu}\\
\ldots \rto& C_4\rto^{\partial_4}&C_3\rto^{\partial_3}& 
L\urto^{\lambda'}\drto_{\lambda}&& P\\
& & & &M\urto_{\mu'}  
\enddiagram
$$
together with actions of $P$ on $L, N, M$ and $C_i$ for $i\geq 3,$ and a function 
$h : M\times N\longrightarrow L.$ The following axioms need to be satisfied.\\
(i) The square $\left(\spreaddiagramrows{-1.2pc} \spreaddiagramcolumns{-1.2pc}
\def\objectstyle{\ssize} \def\labelstyle{\ssize}
\diagram
L \dto_{{\lambda'}} \rto^{\lambda}  & N \dto^{\mu} \\
{M} \rto_{\mu'}& P
\enddiagram\right)$  is a crossed square; \\ 
(ii) The group $C_n$ is abelian for $n \geq 3$ \\
(iii) The boundary homomorphisms satisfy $\partial_n\partial_{n+1} = 1$ for
$n \geq 3,$ and  $\partial_3(C_3)$ lies in the intersection 
$\mbox{ker}\lambda\cap \mbox{ker}{\lambda'};$\\
(iv)  The action of $P$ on $C_n$ for $n \geq 3$ is such that ${\mu}{M}$ and
${\mu'}N$ act trivially. Thus each $C_n$ is a $\pi_0$-module with $\pi_0 =
P/{\mu}{M}{\mu'}N.$ \\
(v) The homomorphisms $\partial_n$ are $\pi_0$-module homomorphisms for $n \geq 3.$

This last condition does make sense  since the axioms for crossed squares imply that 
$\mbox{ker}{\mu'}\cap \mbox{ker}{\mu}$ is a $\pi_0$-module.

A morphism of squared complexes $$\Phi: \bigg(C_\ast,\left(\spreaddiagramrows{-1.2pc} 
\spreaddiagramcolumns{-1.2pc}
\def\objectstyle{\ssize} \def\labelstyle{\ssize}
\diagram
L \dto_{{\lambda'}} \rto^{\lambda}  & N \dto^{\mu} \\
{M} \rto_{\mu'}& P
\enddiagram\right) \bigg)\longrightarrow \bigg({C_\ast'},\left(\spreaddiagramrows{-1.2pc} 
\spreaddiagramcolumns{-1.2pc}
\def\objectstyle{\ssize} \def\labelstyle{\ssize}
\diagram
{L'} \dto_{{\lambda'}} \rto^{\lambda}  & {N'} \dto^{\mu} \\
{M'} \rto_{\mu'}& {P'}
\enddiagram\right)\bigg)$$
consists of a morphism of crossed squares $(\Phi_{L}, 
\Phi_{N}, \Phi_{M}, \Phi_{P})$, together with a family of equivariant homomorphisms $\Phi_n$ for $n \geq 3$ satisfying
$\Phi_{L}\partial_3 = {\partial'}_3\Phi_{L}$ and 
$\Phi_{n-1}\partial_n = {\partial'}_n\Phi_n$ for $n \geq 4.$ 
There is clearly a category $\mathfrak{SqComp}$ of squared complexes. This
exists in both group and groupoid based versions.

A squared complex is thus a crossed square with a `tail' attached.  The same
process can be applied to 2-crossed modules and leads to the notion of a
2-crossed complex, \cite{mp4}.

A 2-crossed complex of group(oid)s is a sequence of group(oid)s 
$$C:\hspace{1cm} \ldots \rightarrow C_n
\stackrel{\partial_n}{\rightarrow}C_{n-1}\rightarrow\quad\ldots \quad C_2
\stackrel{\partial_2}{\rightarrow}C_1\stackrel{\partial_1}{\rightarrow}C_0$$
in which\\
(i) ~$C_n$ is abelian for $n\geq 3$;\\
(ii) ~$C_0$ acts on $C_n$, $n\geq 1$, the action of $\partial C_1$ being
trivial on $C_n$ for $n\geq 3$;\\
(iii)~ each $\partial_n$ is a $C_0$-group(oid) homomorphism and
$\partial_i\partial_{i+1} =1$ for all $i\geq 1$;\\
and\\
(iv)~ $C_2
\stackrel{\partial_2}{\rightarrow}C_1\stackrel{\partial_1}{\rightarrow}C_0$
is a 2-crossed module.

\medskip

The following is an easy consequence of our earlier work.
\begin{thm}
If $$
\diagram
& & & & N\drto^{\mu}\\
\ldots \rto& C_4\rto^{\partial_4}&C_3\rto^{\partial_3}& 
L\urto^{\lambda'}\drto_{\lambda}&& P\\
& & & &M\urto_{\mu'}  
\enddiagram
$$ is a squared complex then $$ \ldots \rightarrow C_n
\stackrel{\partial_n}{\rightarrow}C_{n-1}\rightarrow\ldots C_{3}\rightarrow L
{\rightarrow}N\rtimes M{\rightarrow}P$$is a 2-crossed complex.\hfill$\square$
\end{thm}

\section{Homotopy groups}

It is folklore that any bisimplicial group ${\cal X}_{*,*}$, the diagonal,
$diag ({\cal X})$ and the codiagonal $\nabla {\cal X}$ have the same homotopy
type.  The original work of Loday, \cite{loday} and
Conduch\'e, \cite{conduche}, together with the second author's \cite{porter}
all used the diagonal so showing that the homotopy groups of $diag ({\cal X(M)})$ 
were those of $G$ when ${\cal{M}} = {\mathfrak{M}}(G,2)$.  The use of the mapping
  complex construction by Loday, again in \cite{loday}, also gave the same
  homotopy groups.  Here we will briefly look at the homotopy groups of
  $\nabla({\cal X}({\mathfrak M}(G,2)))$ directly.

The crossed square ${\mathfrak{M}}(G,2)$ as we recalled earlier has form 
$$
\mathfrak{M}(\mathbf{G}, 2) \cong
\left ( \diagram
NG_2 / \partial_3 (NG_3 ) \dto_{\mu_1} \rto^{~~~\mu_2} & \mbox{ Ker }d_{0} 
\dto_{\mu_1} \\
\mbox{ Ker }d_{1} \rto^{~~~\mu_2}  &   G_1    
\enddiagram \right).
$$
The horizontal kernel is $\pi_2(G) \to 1$ and the cokernel ``is'' $\mathfrak{M}(G,1)$, i.e. isomorphic to 
$$\begin{array}{ccc}
NG_1/\partial_2NG_2 \longrightarrow G_0,
\end{array}$$
 which, in turn, has kernel $\pi_1(G)$ and cokernel $\pi_0G)$ .

The corresponding simplicial group $\nabla({\cal X}({\mathfrak M}(G,2)))$ has
Moore complex $$NG_2 / \partial_3 (NG_3 ) \stackrel{(\mu_2^{-1},\mu_1)}{\to}
\mbox{ Ker }d_{0}\rtimes \mbox{ Ker }d_{1} \stackrel{\mu_1\mu_2}{\to}G_1.$$
The homomorphisms $\mu_2$ and $\mu_1$ from $NG_2 / \partial_3 (NG_3 )$ in the
crossed square are both induced from $d_2$, so it is immediate that
$\mbox{Ker}(\mu_2^{-1},\mu_1)$ is $\pi_2(G)$.  The other two $\mu$s (bottom
and 
right of the square) are inclusions.  Since $G_0 \cong
G_1/\mu_2(\mbox{Ker}d_1)$, we again easily check that $\pi_0$ of the complex,
i.e. $G_1/ Im(\mu_1\mu_2)$, is $\pi_0(G)$.

Finally
$$(m,n) \in \mbox{Ker}(\mu_1\mu_2)$$if and only if $\mu_1(m) = \mu_2(n)^{-1}$, 
but as $\mu_1$ and $\mu_2$ are inclusions, this amounts to
$\mbox{Ker}(\mu_1\mu_2)$ being
$$\{(m,m^{-1}) : m \in \mbox{Ker} d_0 \cap \mbox{Ker} d_1\}.$$Then the image of
$(\mu_2^{-1},\mu_1)$ identifies as $$\{(m,m^{-1}) : m \in \partial_2(NG_2)\}$$ 
and again a routine calculation shows that
$$\frac{\mbox{Ker}(\mu_1\mu_2)}{\mbox{Im}(\mu_2^{-1},\mu_1)} \cong
\pi_1(G),$$ as expected.

These calculations do not involve the Peiffer lifting so are not conclusive
about homotopy type, however a neat argument noted by Conduch\'e shows that
there is an (obvious) epimorphism
$$\xymatrix{NG_2/\partial_3(NG_3)\ar[d]_= \ar[r]&\mbox{ Ker }d_{0}\rtimes \mbox{ Ker }d_{1}\ar[d] \ar[r]&G_1\ar[d]\\
NG_2/\partial_3(NG_3)\ar[r]&\mbox{ Ker }d_{0}\ar[r] &G_0}$$
with acyclic kernel.  This still does not quite clinch the argument, since it would
be better to check that there was an acyclic fibration
$$\nabla({\cal X}({\mathfrak M}(G,2)))\to t_{2]}G$$where $t_{2]}G$ here
denotes the homotopy truncation of $G$ (essentially its 2-coskeleton).  However 
we have not attempted to give this here.

\section{Higher dimensions}

Loday's mapping complex was defined for cat$^n$-groups (see \cite{loday}) and
results of Bullejos, Cegarra and Duskin, \cite{manantjack} suggest that a
similar multiple codiagonal would give his mapping complex.  This raises the
question of what would be the result on  taking a crossed $n$-cube $\cal M$ and forming its $n$-fold nerve 
$Ner^{(n)}{\cal{M}}$, which will be an $n$-simplicial group (just generalise
the construction of section 6.1).  Again Bullejos, Cegarra and Duskin use an
inductive argument to derive a result that would suggest the multicodiagonal
of $Ner^{(n)}{\cal{M}}$ should have Moore complex of length $n$.

We first note that the proof of Proposition 6.1 shows the following:

\begin{prop}

If ${\cal X}_{*,*}$ is a bisimplicial group such that for any $p$, $N({\cal
  X}_{p,*})_q = 1$ for $q \geq 2$, whilst for any $q$, $N({\cal
  X}_{*,q})_p = 1$ for $p \geq m$, then 
$$N(\nabla{\cal X})_n = 1 \mbox{~~~ for } n \geq m + 1.$$ \hfill$\blacksquare$
\end{prop}

We next need a categorical description of the multidiagonal:

Let $ or : \Delta \times \Delta \to \Delta$ be the ordinal sum functor, then
for a bisimplicial group ${\cal X}$, it is well known that $\nabla {\cal X}$
has a description as a coend
$$(\nabla {\cal X})_n = \int^{[p],[q]}\Delta([n],[p]or[q]) \times {\cal
  X}_{p,q},$$
(cf. for example, Cordier-Porter, \cite{JMCTP}). A corresponding codiagonal
for a $m$-fold simplicial group ${\cal X}_{\underline{*}}$, $\underline{*}$ ~ an
$m$-fold index, is  
$$(\nabla^{(m)} {\cal X})_n = \int^{\underline{p}}\Delta([n],or\underline{p})
\times {\cal X}_{\underline{p}},$$
where by abuse of notation, we indicate by $$or : \Delta^{\times m} \to
\Delta,$$
the $m$-fold ordinal sum, $\underline{p} = ([p_1], \ldots, [p_m])$ an $m$-fold 
index and have used as a shorthand $$or \underline{p} = [p_1]or \ldots or
[p_m],$$which causes no ambiguity as $or$ is associative.

As $or \underline{p} = ([p_1]or \ldots or [p_{m-1}]) or [p_m] =: (or
\underline{p^\prime}) or[p_m]$
$$(\nabla^{(m)} {\cal X})_n = \int^{\underline{p^\prime}}\int^{[p_m]}\Delta([n],(or\underline{p^\prime})or[p_m])
\times {\cal X}_{\underline{p^\prime,p_m}},$$
but then `integrating' over all $\underline{p^\prime}$ for each $n$ produces a 
new description of $(\nabla^{(m)} {\cal X})$ as $\nabla ( \nabla^{(m-1)} {\cal 
  X}_{\underline{*}, ?})$ i.e. an iterative description.  (The proof uses the
fact that any mapping from $[n]$ to $(or\underline{p^\prime})or[p_m]$
effectively partitions $[n]$ into an initial segment mapping to
$or\underline{p^\prime}$ and a second part mapping to $[p_m]$. Fixing the
latter part, we form the coend with $[n]$ replaced by the first segment of the 
partition.)  Now assume that the $m$-simplicial group ${\cal
  X}_{\underline{*}}$ is obtained as the $m$-fold nerve of a crossed $m$-cube
(or cat$^m$-group), then applying the above proposition repeatedly we obtain:

\begin{thm}
If $\mathfrak{M}$ is a crossed $m$-cube with $m$-fold nerve, the
$m$-simplicial group, ${\cal X}(\mathfrak{M})$, then 
$$N(\nabla{\cal X}(\mathfrak{M}))_n = 1 \quad \mbox {for} \quad n \geq
m+1$$\hfill $\blacksquare$
\end{thm}

There is a `complex' form of this result as well.  If one takes the obvious
notion of $m$-cube complex, generalising the squared complexes considered
above, then Theorem 6.2 generalises without bother to give a notion of
$m$-crossed complex and a construction that generalises $\nabla {\cal X}$
above.  The Moore complex that one gets is related to complexes considered by
Duskin and Nan Tie as well as the hypercrossed complexes of \cite{C&C1}.

The structural maps of such hypercrossed complexes and the related Peiffer
pairings considered in earlier papers of the series could be given explicit
algebraic formulae in terms of the $h$-maps in the $m$-cube complex, but the
problem of the complexity of these formulae raises doubts as to their usefulness.

\small{
{\noindent  A. Mutlu}\hspace{5.55cm} \ \ \ T. Porter  \\      
{\it Department of Mathematics \ \ \hspace{2.5cm}School of Informatics\\
Faculty of Science \hspace{3.6cm} \ \ \ \ \ University of Wales Bangor, \\
University of Celal Bayar \hspace{3cm}\ Gwynedd, LL57 1UT, UK.\\
Manisa, Turkey \hspace{4.75cm} e-Mail: t.porter@bangor.ac.uk\\
e-Mail: amutlu@spil.bayar.edu.tr\\}

\end{document}